# A POLYNOMIAL TIME ALGORITHM FOR VERTEX ENUMERATION AND OPTIMIZATION OVER SHAPED PARTITION POLYTOPES

FRANK K. HWANG, SHMUEL ONN, AND URIEL G. ROTHBLUM

ABSTRACT. We consider the *Shaped Partition Problem* of partitioning $n$ given vectors in real $k$-space into $p$ parts so as to maximize an arbitrary objective function which is convex on the sum of vectors in each part, subject to arbitrary constraints on the number of elements in each part. In addressing this problem, we study the *Shaped Partition Polytope* defined as the convex hull of solutions. The Shaped Partition Problem captures $\mathcal{NP}$-hard problems such as the Max-Cut problem and the Traveling Salesperson problem, and the Shaped Partition Polytope may have exponentially many vertices and facets, even when $k$ or $p$ are fixed. In contrast, we show that when both $k$ and $p$ are fixed, the number of vertices is polynomial in $n$, and all vertices can be enumerated and the optimization problem solved in strongly polynomial time. Explicitly, we show that any Shaped Partition Polytope has $O(n^{k\binom{p}{2}})$ vertices which can be enumerated in $O(n^{k^2 p^3})$ arithmetic operations, and that any Shaped Partition Problem is solvable in $O(n^{kp^2})$ arithmetic operations.

## 1. INTRODUCTION

The *Partition Problem* concerns the partitioning of vectors $A^1, \ldots, A^n$ in $k$-space into $p$ parts so as to maximize an objective function which is convex on the sum of vectors in each part. Each vector $A^i$ represents $k$ numerical attributes associated with the $i$th element of the set $[n] = \{1, \ldots, n\}$ to be partitioned. Each ordered partition $\pi = (\pi_1, \ldots, \pi_p)$ of $[n]$ is then associated with the $k \times p$ matrix $A^\pi = \left[\sum_{i \in \pi_1} A^i, \ldots, \sum_{i \in \pi_p} A^i\right]$ whose $j$th column represents the total attribute vector of the $j$th part. The problem is to find an *admissible* partition $\pi$ which maximizes an objective function $f$ given by $f(\pi) = C(A^\pi)$ where $C$ is a real convex functional on $\mathbb{R}^{k \times p}$. Of particular interest is the *Shaped Partition Problem*, where the admissible partitions are those $\pi$ whose *shape* $(|\pi_1|, \ldots, |\pi_p|)$ lies in a prescribed set $\Lambda$ of admissible shapes. In this article we concentrate on this later situation.

The Shaped Partition Problem has applications in diverse fields that include circuit layout, clustering, inventory, splitting, ranking, scheduling and reliability, see

The research of Shmuel Onn was supported by the Mathematical Sciences Research Institute at Berkeley California through NSF Grant DMS-9022140, and by Two Technion VPR Grants.

The research of Uriel G. Rothblum was supported by the E. and J. Bishop Research Fund at the Technion and by ONR Grant N00014-92-J1142.

[3, 4, 6, 11] and references therein. As we demonstrate later on, the problem has expressive power that captures $\mathcal{NP}$-hard problems such as the Max-Cut problem and the Traveling Salesperson problem, even when the number $p$ of parts or attribute dimension $k$ are fixed.

Our first goal in this article is to show that a polynomial time algorithm for the Shaped Partition Problem *does* exist when both $p$ and $k$ are fixed. This result remains valid when the set $\Lambda$ of admissible shapes and the function $C$ are presented by oracles (see Section 4 for the precise statement):

- **Theorem 4.6:** Any Shaped Partition Problem is solvable in oracle time $O(n^{kp^2})$.

Since $C$ is convex, the Shaped Partition Problem can be embedded into the problem of maximizing $C$ over the *Shaped Partition Polytope* $\mathcal{P}_A^\Lambda$ defined to be the convex hull of all matrices $A^\pi$ corresponding to partitions of admissible shapes. The class of Shaped Partition Polytopes is very broad and generalizes and unifies classical permutation polytopes such as Birkhoff's polytope and the Permutohedron (see e.g. [15]). Its subclass of *bounded* shape partition polytopes with lower and upper bounds on the shapes was previously studied in [1], under the assumption that the vectors $A^1, \ldots, A^n$ are *distinct*. Therein a procedure for testing whether a given $A^\pi$ is a vertex of $\mathcal{P}_A^\Lambda$ was obtained. This procedure is simplified and extended in [10]. A related but different important generalization of classical permutation polytopes, arizing when algebraic (representation-theoretical) constraints, rather than shape constraints, are enforced on the permutations involved, was studied in [15] and references therein.

Since a Shaped Partition Polytope is defined as the convex hull of an implicitly presented set whose size is typically exponential in the input size even when both $p$ and $k$ are fixed, neither a presentation as the convex hull of vertices nor as the intersection of halfspaces is readily available for it. Our second objective is to bypass this difficulty by proving that, nevertheless, for fixed $p$ and $k$, the number of vertices of Shaped Partition Polytopes *is* polynomially bounded in $n$, and it is possible to explicitly enumerate all vertices in polynomial time (see Section 4 for the precise statements):

- **Theorem 4.3:** Any Shaped Partition Polytope $\mathcal{P}_A^\Lambda$ has $O(n^{(k+1)\binom{p}{2}})$ vertices.

While this bound suffices for our algorithmic considerations in this article, we point out that in [9] we establish the more accurate bound of $O(n^{k\binom{p}{2}})$ on the number of vertices.

- **Theorem 4.5:** All vertices of $\mathcal{P}_A^\Lambda$ can be enumerated in oracle time $O(n^{k^2p^3})$.

In particular, it is possible to compute the *number* of vertices efficiently. This might be extendible to the situation of variable $k$ and $p$, while counting vertices is generally a hard task (see [12, 13]), as is counting partitions with various prescribed properties (see [2]). The vertex counting problem for variable $k$ and $p$ will be addressed elsewhere.

The article is organized as follows. In the next section we formally define the Shaped Partition Problem and Shaped Partition Polytope. We demonstrate that the Max-Cut problem and Traveling Salesperson problem can be modeled as Shaped Partition Problems with fixed $p = 2$ and $k=1$, respectively, and that the corresponding polytopes have exponentially many vertices. In Section 3 we establish a sequence of lemmas in preparation for Section 4. We extend earlier results of [1] and provide a necessary condition for $A^\pi$ to be a vertex of $\mathcal{P}_A^\Lambda$. This leads to the study of partitions of points in space whose convex hulls are pairwise separable (in the special case $k = p = 2$ such partitions had been studied quite extensively, see e.g. [14]; the case $k = 3, p = 2$ has also been considered quite recently in [5]). In the final Section 4 we reward our preparatory results of Section 3 and establish Theorems 4.3, 4.5, and 4.6.

## 2. Shaped Partition Problems and Shaped Partition Polytopes

For a matrix $A$, we use $A^i$ to denote its $i$th column. A *p-partition* of $[n] := \{1, \ldots, n\}$ is an ordered collection $\pi = (\pi_1, \ldots, \pi_p)$ of $p$ (possibly empty) sets whose disjoint union is $[n]$. A *p-shape* of $n$ is a tuple $\lambda = (\lambda_1 \ldots, \lambda_p)$ of nonnegative integers $\lambda_1 \ldots, \lambda_p$ satisfying $\sum_{i=1}^p \lambda_i = n$. The *shape of a partition* $\pi$ is $|\pi| := (|\pi_1|, \ldots, |\pi_p|)$. If $\Lambda$ is a set of $p$-shapes of $n$ then a $\Lambda$-*partition* is any partition $\pi$ whose shape $|\pi|$ is a member of $\Lambda$.

Let $A$ be a real $k \times n$ matrix. For each $p$-partition $\pi$ of $[n]$ we define the *A-matrix* of $\pi$ to be the $k \times p$ matrix

$$A^\pi = \left[ \sum_{i \in \pi_1} A^i, \ldots, \sum_{i \in \pi_p} A^i \right],$$

where $\sum_{i \in \pi_j} A^i := 0$ whenever $\pi_j = \emptyset$. We consider the following algorithmic problem.

**Shaped Partition Problem.** Given positive integers $p, k, n$, matrix $A \in \mathbb{R}^{k \times n}$, nonempty set $\Lambda$ of $p$-shapes of $n$, and objective function on $\Lambda$-partitions given by $f(\pi) = C(A^\pi)$ with $C$ convex on $\mathbb{R}^{k \times p}$, find a $\Lambda$-partition $\pi^*$ which maximizes $f$, namely satisfies

$$f(\pi^*) = \max\{f(\pi) : \; |\pi| \in \Lambda\}.$$

The complexity of this problem depends on the presentation of $\Lambda$ and $C$ and on whether real arithmetic computations or Turing machine computations on rational data are assumed. Pleasantly, our algorithms work in strongly polynomial time and can cope with minimal information on $\Lambda$ and $C$. So both $\Lambda$ and $C$ may be presented by oracles. All that is needed of an oracle for $\Lambda$ is that on query shape $\lambda$ it answers whether or not $\lambda \in \Lambda$. All that is needed of an oracle for $C$ is that on query $A^\pi$ it answers with $C(A^\pi)$.

Since $C$ is convex, the Shaped Partition Problem can be embedded into the problem of maximizing $C$ over the convex hull of $A$-matrices of feasible partitions, defined as follows.

**Shaped Partition Polytope.** For a matrix $A \in \mathbb{R}^{k \times n}$ and nonempty set $\Lambda$ of $p$-shapes of $n$ we define the *Shaped Partition Polytope* $\mathcal{P}_A^\Lambda$ to be the convex hull of all $A$-matrices of $\Lambda$-partitions

$$\mathcal{P}_A^\Lambda := \mathrm{conv}\,\{A^\pi : \ |\pi| \in \Lambda\} \subset \mathbb{R}^{k \times p}.$$

We point out that $\mathcal{P}_A^\Lambda$ is the image of the Shaped Partition Polytope $P_I^\Lambda$, with $I$ the $n \times n$ identity, under the projection $X \mapsto A \cdot X$. In [8] this is exploited, for the situation where the function $C$ is *linear* and $\Lambda = \{\lambda : \ l \leq \lambda \leq u\}$ is a set of *bounded* shapes, to solve the corresponding Shaped Partition Problem for all $n, p, k$ in polynomial time by linear programming over $P_I^\Lambda$, and to efficiently test for a given partition $\pi$ whether $A^\pi$ is a vertex of $\mathcal{P}_A^\Lambda$.

Next we demonstrate that, even when either $k$ or $p$ are fixed, the Shaped Partition Problem may be $\mathcal{NP}$-hard, and the number of vertices of the Shaped Partition Polytope may be exponential. Therefore, polynomial time algorithms for optimization and vertex enumeration may (and, as we show, *do*) exist only when *both* $k$ and $p$ are fixed.

**Example 2.1. Max-Cut Problem and Unit Cube.** Fix $p = 2$. Let $\Lambda$ be the set of all 2-shapes of $n$, let $k = n$ and let $A$ be the $n \times n$ identity matrix. The $A$-matrices of $\Lambda$-partitions in this case are precisely all $(0, 1)$-valued $n \times 2$ matrices with each row sum equals 1. The Shaped Partition Polytope $\mathcal{P}_A^\Lambda$ has $2^n$ vertices which stand in bijection with $\Lambda$-partitions, and is affinely equivalent to the $n$-dimensional unit cube by projection of matrices onto their first column. Now let $G = ([n], E)$ be an undirected graph, and define an objective function $f_G$ as follows: for a 2-partition $\pi = (\pi_1, \pi_2)$ of $[n]$ let $f_G(\pi) = |\{e \in E : \ |e \cap \pi_1| = 1\}|$ be the number of edges in the *cut* of $G$ defined by $\pi$. Since each $A^\pi$ is a distinct vertex of $\mathcal{P}_A^\Lambda$, there is a convex $C$ on $\mathbb{R}^{k \times 2}$ such that $f_G(\pi) = C(A^\pi)$ for all $\pi$. The Shaped Partition Problem in this case is precisely the Max-Cut problem of finding a largest cut in a graph hence is $\mathcal{NP}$-hard.

**Example 2.2. Travelling Salesperson Problem and Permutohedron.** Fix $k = 1$. Let $p = n$, let $\Lambda = \{1^n\}$ consist of the single shape $1^n = (1, \ldots, 1)$, and let $A = [1, 2, \ldots, n]$. The matrices $A^\pi$ in this case are simply all permutations of $A$. The Shaped Partition Polytope $\mathcal{P}_A^\Lambda$ has $n!$ vertices which stand in bijection with $\Lambda$-partitions, and is the so-called *Permutohedron*. Now let $D \in \mathbb{R}^{n \times n}$ be a symmetric nonnegative matrix, where $D_{i,j}$ represents the *distance* between locations $i$ and $j$.

Define an objective function $f_D$ as follows. For partition $\pi = (\{i_1\}, \ldots, \{i_n\})$ let $f_D(\pi) = -\sum_{j=1}^n D_{i_j, i_{j+1}}$ be minus the length of the *Traveling Salesperson Tour* on $[n]$ given by $\pi$ (with $i_{n+1} := i_1$). Since each $A^\pi$ is a distinct vertex of $\mathcal{P}_A^\Lambda$, there is a convex $C$ on $\mathbb{R}^n$ such that $f_D(\pi) = C(A^\pi)$ for all $\pi$. The Shaped Partition Problem in this case is precisely the Traveling Salesperson Problem of finding a shortest tour hence is $\mathcal{NP}$-hard.

We finish this section with an example demonstrating how our results enable a polynomial time solution procedure which bypasses the need for checking the exponentially many feasible solutions present.

**Example 2.3. Splitting.** The $n$ assets of a company are to be split among its $p$ owners as follows. For $i = 1, \ldots, p$, the $i$th owner prescribes a nonnegative vector $A_i = (A_{i,1}, \ldots, A_{i,n})$ with $\sum_{j=1}^n A_{i,j} = 1$, whose entries represent the relative values of the various assets to this owner. A partition $\pi = (\pi_1, \ldots, \pi_p)$ is sought which splits the assets among the owners as evenly as possible, and maximizes the $l_q$-norm $(\sum_{i=1}^p |\sum_{j \in \pi_i} A_{i,j}|^q)^{\frac{1}{q}}$ of the vector whose $i$th entry $\sum_{j \in \pi_i} A_{i,j}$ is the total relative value of the assets handed to the $i$th owner by $\pi$. This can be formulated as a Shaped Partition Problem with $k = p$, the $p \times n$ matrix $A$ whose rows are the $A_i$, the set $\Lambda$ of all $p$-shapes $\lambda$ of $n$ with $\lfloor \frac{n}{p} \rfloor \leq \lambda_i \leq \lceil \frac{n}{p} \rceil$ for all $i$, and objective function $f_q(\pi) = C_q(A^\pi)$ where $C_q$ is the convex function on $\mathbb{R}^{p \times p}$ defined by $C_q(M) = \sum_{i=1}^p |M_{i,i}|^q$. For fixed $p$, by Theorem 4.6 we can find an optimal partition in polynomial time $O(n^{p^3})$, while the number of partitions is exponential $\Omega(p^n n^{\frac{1-p}{2}})$. We note further that for $q = 1$ the functional $C_1$ is linear on $\mathbb{R}_+^{p \times p}$ hence our results of [8] apply and yield a polynomial time solution even when $p$ is variable.

## 3. Preparation

In this section we establish a sequence of lemmas in preparation for the final section. Two finite sets $U, V$ of points in $\mathbb{R}^k$ will be called *separable* if there is a vector $h \in \mathbb{R}^k$ such that $h^T u < h^T v$ for all $u \in U$ and $v \in V$ with $u \neq v$. It is not hard to see that $U, V$ are separable if and only if their convex hulls $\text{conv}(U)$ and $\text{conv}(V)$ are either disjoint or intersect in a single point which is a common vertex of both. In particular, if $U$ and $V$ are separable then $|U \cap V| \leq 1$.

Let $A$ be a given $k \times n$ matrix. For a subset $S \subseteq [n]$ let $A^S = \{A^i : i \in S\}$ be the set of columns of $A$ indexed by $S$ (with multiple copies of columns identified). A $p$-partition $\pi = (\pi_1, \ldots, \pi_p)$ will be called $A$-*separable* if the sets $A^{\pi_r}$ and $A^{\pi_s}$ are separable for each pair $1 \leq r < s \leq p$. So $\pi$ is $A$-separable if for each pair $1 \leq r < s \leq p$ there is a vector $h_{r,s} \in \mathbb{R}^k$ such that $h_{r,s}^T A^i < h_{r,s}^T A^j$ for all $i \in \pi_r$ and $j \in \pi_s$ with $A^i \neq A^j$. We have the following lemma which, for a matrix $A$ with *distinct nonzero* columns, specializes to a result of [1].

**Lemma 3.1.** *Let $A$ be a matrix in $\mathbb{R}^{k \times n}$, let $\Lambda$ be a nonempty set of $p$-shapes of $n$, and let $\pi$ be a $\Lambda$-partition. If $A^\pi$ is a vertex of $\mathcal{P}_A^\Lambda$ then $\pi$ is an $A$-separable partition.*

*Proof.* The claim being obvious for $p = 1$, suppose that $p \geq 2$. Let $A^\pi$ be a vertex of $\mathcal{P}_A^\Lambda$. Then there is a matrix $C \in \mathbb{R}^{k \times p}$ such that the linear functional on $\mathbb{R}^{k \times p}$ given by the inner product $\langle C, X \rangle = \sum_{i=1}^k \sum_{j=1}^p C_i^j X_i^j$ is uniquely maximized over $\mathcal{P}_A^\Lambda$ at $A^\pi$. Pick any pair $1 \leq r < s \leq p$, and let $h_{r,s} = C^s - C^r$. Suppose there are $i \in \pi_r$ and $j \in \pi_s$ with $A^i \neq A^j$ (otherwise $A^{\pi_r}$ and $A^{\pi_s}$ are trivially separable). Let $\pi'$ be the $\Lambda$-partition obtained from $\pi$ by switching $i$ and $j$, i.e. taking $\pi'_r := \pi_r \cup \{j\} \setminus \{i\}$, $\pi'_s := \pi_s \cup \{i\} \setminus \{j\}$, and $\pi'_t := \pi_t$ for all $t \neq r, s$. Then $(A^{\pi'})^r = (A^\pi)^r + A^j - A^i \neq (A^\pi)^r$ hence $A^{\pi'} \neq A^\pi$. By choice of $C$ we have $\langle C, A^{\pi'} \rangle < \langle C, A^\pi \rangle$ and so

$$h_{r,s}^T(A^j - A^i) = (C^s - C^r)^T(A^j - A^i) = \sum_{t=1}^p (C^t)^T((A^\pi)^t - (A^{\pi'})^t) = \langle C, A^\pi - A^{\pi'} \rangle > 0.$$

This proves that $A^{\pi_r}$ and $A^{\pi_s}$ are separable for each pair $1 \leq r < s \leq p$, hence $\pi$ is $A$-separable. □

Let $A$ be a $k \times n$ matrix. A $p$-partition $\pi = (\pi_1, \ldots, \pi_p)$ will be called $A$-*disjoint* if $\mathrm{conv}(A^{\pi_r})$ and $\mathrm{conv}(A^{\pi_s})$ are disjoint for each pair $1 \leq r < s \leq p$. Also, for $A \in \mathbb{R}^{k \times n}$ let $\underline{A}$ be the matrix obtained by appending a new row $[1, \ldots, n]$ to $A$.

**Lemma 3.2.** *Let $A \in \mathbb{R}^{k \times n}$ and let $\underline{A}$ be as above. For each $A$-separable partition $\pi$ there exists an $\underline{A}$-disjoint partition $\underline{\pi}$ of the same shape $|\underline{\pi}| = |\pi|$ satisfying $A^\pi = A^{\underline{\pi}}$.*

*Proof.* Let $\pi = (\pi_1, \ldots, \pi_p)$ be any $A$-separable $p$-partition. As any 1-partition is both separable and disjoint, we may assume $p \geq 2$. Let $\underline{\pi}$ be the *lexicographically first* partition for which $\{A^{\underline{\pi}_i}\} = \{A^{\pi_i}\}$ for $i = 1, \ldots, p$, where $\{A^S\}$ stands for the *multiset* of columns of $A$ indexed by $S$. More precisely, $\underline{\pi}$ is the partition such that for all $u \in \mathbb{R}^k$ the following hold: if we put $U_i := \{j \in \pi_i : A^j = u\}$ and $\underline{U}_i := \{j \in \underline{\pi}_i : A^j = u\}$, then $|\underline{U}_i| = |U_i|$ for all $i = 1, \ldots, p$, and $i < j$ for all $i \in \underline{U}_r$ and $j \in \underline{U}_s$ with $1 \leq r < s \leq p$. By definition $\underline{\pi}$ has the same shape as $\pi$ and satisfies $A^{\underline{\pi}} = A^\pi$.

Now consider any pair $1 \leq r < s \leq p$. Since $\pi$ is $A$-separable, there is a vector $h_{r,s} \in \mathbb{R}^k$ such that $h_{r,s}^T A^i < h_{r,s}^T A^j$ for all $i \in \pi_r$ and $j \in \pi_s$ with $A^i \neq A^j$. Let $\underline{h}_{r,s}$ be the $(k+1)$-vector obtained by appending a last coordinate $\epsilon$ to $h_{r,s}$, with $\epsilon > 0$ sufficiently small. For $i \in \underline{\pi}_r$ and $j \in \underline{\pi}_s$ consider the quantity

$$\underline{h}_{r,s}^T(\underline{A}^j - \underline{A}^i) = h_{r,s}^T(A^j - A^i) + \epsilon(j - i);$$

if $A^i \neq A^j$ then, since $\epsilon$ is infinitesimal, this quantity is positive; if $A^i = A^j$ then, by choice of $\underline{\pi}$, $i < j$ and this quantity is again positive. So $\mathrm{conv}(\underline{A}^{\underline{\pi}_r})$ and $\mathrm{conv}(\underline{A}^{\underline{\pi}_s})$ are disjoint. It follows that $\underline{\pi}$ is $\underline{A}$-disjoint and the proof is complete. □

For $v \in \mathbb{R}^k$ denote by $\bar{v}$ the $(k+1)$-vector obtained by appending a first coordinate 1 to $v$. An *orientation* of a hyperplane $H$ is a presentation of $H$ as the zero set $H = \{x \in \mathbb{R}^k : h^T \bar{x} = 0\}$ of a linear form (with $h \in \mathbb{R}^{k+1}$). With an orientation of $H$, the closed and open halfspaces $H^\leq, H^<$ *below* it are distinguishable from the ones $H^\geq, H^>$ *above* it. The hyperplane $H$ spanned by an *ordered* list $[v^1, \ldots, v^k]$ of $k$ affinely independent points in $\mathbb{R}^k$ will be oriented as $H = \{x \in \mathbb{R}^k : \det[\bar{v}^1, \ldots, \bar{v}^k, \bar{x}] = 0\}$.

A matrix $A$ will be called *generic* if its columns are in affine general position. Let $A \in \mathbb{R}^{k \times n}$ be generic with $n > k$. With each triple $(I, J^<, J^>)$ where $I$ is a $k$-subset of $[n]$ and $(J^<, J^>)$ is a 2-partition of $I$, we associate two $A$-disjoint 2-partitions $\pi^<$ and $\pi^>$ as follows. Let $[A^I]$ denote the list of columns of $A$ indexed by $I$ and ordered as in $A$. Let $H := \mathrm{aff}([A^I])$ be the hyperplane spanned by $A^I$ oriented by the ordered list $[A^I]$. Let $I^< := \{i \in [n] : A^i \in H^<\}$ and $I^> := \{i \in [n] : A^i \in H^>\}$. Define $\pi^< := (I^< \cup J^<, I^> \cup J^>)$ and $\pi^> := (I^> \cup J^>, I^< \cup J^<)$. Both $\pi^<$ and $\pi^>$ are $A$-disjoint since $H$ can be slightly perturbed so as to strictly separate $A^{\pi_1}$ and $A^{\pi_2}$.

**Lemma 3.3.** *Let $A \in \mathbb{R}^{k \times n}$ be generic. The set of $A$-disjoint 2-partitions is the set of all 2-partitions if $n \leq k$ and is the set of partitions associated with triples $(I, J^<, J^>)$ as above if $n > k$.*

*Proof.* If $n \leq k$ then the columns of $A$ are affinely independent hence $\mathrm{conv}(A^{\pi_1})$ and $\mathrm{conv}(A^{\pi_2})$ are disjoint for any 2-partition $\pi$. Suppose now $n > k$ and let $\pi = (\pi_1, \pi_2)$ be any $A$-disjoint 2-partition. Since $\mathrm{conv}(A^{\pi_1})$ and $\mathrm{conv}(A^{\pi_2})$ are disjoint, a standard separation result assures the existence of a hyperplane spanned by columns of $A$ such that $A^{\pi_1}$ is contained in one of its closed halfspaces and $A^{\pi_2}$ in the other. Since $A$ is generic, this hyperplane $H$ contains exactly $k$ columns of $A$, so that $H = \mathrm{aff}(A^I)$ for some $k$-subset $I$ of $[n]$. If the orientation $\mathrm{aff}([A^I])$ of $H$ is such that $A^{\pi_1} \subset H^\leq$ and $A^{\pi_2} \subset H^\geq$ then $\pi$ is the partition $\pi^<$ associated with the triple $(I, I \cap \pi_1, I \cap \pi_2)$, whereas if $A^{\pi_2} \subset H^\leq$ and $A^{\pi_1} \subset H^\geq$ then $\pi$ is the partition $\pi^>$ associated with the triple $(I, I \cap \pi_2, I \cap \pi_1)$. □

Let now $A \in \mathbb{R}^{k \times n}$ be an arbitrary matrix. Define for each $\epsilon > 0$ a perturbation $A(\epsilon) \in \mathbb{R}^{k \times n}$ of $A$ as follows: for $i = 1, \ldots, n$ let the $i$th column of $A(\epsilon)$ be $A(\epsilon)^i := A^i + \epsilon M_k^i$ where $M_k^i := [i, i^2, \ldots, i^k]^T$ is the image of $i$ on the moment curve in $\mathbb{R}^k$.

**Definition 3.4.** Let $A \in \mathbb{R}^{k \times n}$ and $p \geq 1$. A $p$-partition $\pi$ of $[n]$ will be called $A$-*generic* if it is $A(\epsilon)$-disjoint for all sufficiently small $\epsilon > 0$. The set of $A$-generic $p$-partitions will be denoted by $\Pi_A^p$.

**Lemma 3.5.** *Let $A \in \mathbb{R}^{k \times n}$ and $p \geq 1$. For all sufficiently small $\epsilon > 0$, the perturbation $A(\epsilon)$ is generic and the set of $A(\epsilon)$-disjoint $p$-partitions contains the set of $A$-disjoint $p$-partitions and is independent of $\epsilon$ and equals $\Pi_A^p$.*

*Proof.* Regard the entries of $A(\epsilon)$ as linear forms in $\epsilon$. Recall that $\bar{v}$ denotes the $(k+1)$-vector obtained by appending a first coordinate 1 to $v \in \mathbb{R}^k$. Suppose first that $n > k$. Then each determinant $\det[\bar{A}(\epsilon)^{i_0}, \ldots, \bar{A}(\epsilon)^{i_k}]$ is a polynomial of degree $k$ in $\epsilon$ whose leading coefficient is the nonzero Van der Monde determinant $\det[\bar{M}_k^{i_0}, \ldots, \bar{M}_k^{i_k}]$, hence has finitely many zeroes. So for all sufficiently small $\epsilon > 0$ all of these determinants are nonzero hence $A(\epsilon)$ is generic. For $n \leq k$ the same argument applies by restricting attention to the first $n$ rows of $A$.

Next consider any $A$-disjoint $p$-partition $\pi$. The claim being obvious for $p = 1$, suppose that $p \geq 2$. For each pair $1 \leq r < s \leq p$, there is a hyperplane which strictly separates $A^{\pi_r}$ and $A^{\pi_s}$. Clearly, for all sufficiently small $\epsilon$ this hyperplane also strictly separates $A(\epsilon)^{\pi_r}$ and $A(\epsilon)^{\pi_s}$, so $\pi$ is also $A(\epsilon)$-disjoint.

Finally, we show that for all sufficiently small $\epsilon > 0$ the set of $A(\epsilon)$-disjoint $p$-partitions is the same. We start with the case $p = 2$. If $n \leq k$ then by Lemma 3.3 this set simply consists of all 2-partitions. So assume that $n > k$. Let $I = \{i_1, \ldots, i_k\}$ be any $k$-subset of $[n]$ with $i_1 < \cdots < i_k$. For each $i \notin I$, consider the following determinant, which is again a nonzero polynomial of degree $k$ in $\epsilon$:

$$D^i(\epsilon) := \det[\bar{A}(\epsilon)^{i_1}, \ldots, \bar{A}(\epsilon)^{i_k}, \bar{A}(\epsilon)^i] = \sum_{j=0}^{k} D_j^i \epsilon^j.$$

For all sufficiently small $\epsilon > 0$ the sign of $D^i(\epsilon)$ equals the sign of the first nonzero coefficient among $D_0^i, \ldots, D_k^i$ and is either positive or negative and independent of $\epsilon$. Denote this *generic* sign by $\text{sign}(D^i)$. Then $I^< = \{i \in [n] : \text{sign}(D^i) < 0\}$ and $I^> = \{i \in [n] : \text{sign}(D^i) > 0\}$ are independent of $\epsilon$ and therefore for each 2-partition $(J^<, J^>)$ of $I$, the partitions $\pi^<, \pi^>$ associated with the triple $(I, J^<, J^>)$ are independent of $\epsilon$. So the set of partitions associated with triples $(I, J^<, J^>)$ is the same for all small $\epsilon > 0$ and by Lemma 3.3 this is the set of $A(\epsilon)$-disjoint 2-partitions. Now consider any $p$-partition $\pi$ with $p \geq 2$ arbitrary. Consider any pair $1 \leq r < s \leq p$, and restrict attention to the submatrix $[A^{\pi_r \cup \pi_s}]$ of $A$ consisting of the columns of $A$ indexed by $\pi_r \cup \pi_s$. The same determinantal argument as in the paragraph above shows that for all sufficiently small $\epsilon$ the 2-partition $(\pi_r, \pi_s)$ of $\pi_r \cup \pi_s$ is either $[A(\epsilon)^{\pi_r \cup \pi_s}]$-disjoint or not independently of $\epsilon$. So $\text{conv}(A(\epsilon)^{\pi_r})$ and $\text{conv}(A(\epsilon)^{\pi_s})$ are either disjoint or not independently of $\epsilon$. Since this hold for each pair $1 \leq r < s \leq p$, we see that $\pi$ is either $A(\epsilon)$-disjoint or not independently of $\epsilon$ as claimed. $\square$

Let $p \geq 2$. With each list $[\pi^{r,s} = (\pi_1^{r,s}, \pi_2^{r,s}) : 1 \leq r < s \leq p]$ of $\binom{p}{2}$ 2-partitions of $[n]$ associate a $p$-tuple $\pi = (\pi_1, \ldots, \pi_p)$ of subsets of $[n]$ as follows: for $i = 1, \ldots, p$ put

$$\pi_i := \left(\cap_{j=i+1}^{p} \pi_1^{i,j}\right) \bigcap \left(\cap_{j=1}^{i-1} \pi_2^{j,i}\right).$$

Since $\pi_r \subseteq \pi_1^{r,s}$ and $\pi_s \subseteq \pi_2^{r,s}$ for all $1 \leq r < s \leq p$, the $\pi_i$ are pairwise disjoint. If also $\cup_{i=1}^{p} \pi_i = [n]$ then $\pi$ is a partition, and will be called the *partition* associated with the given list. If moreover each $\pi^{r,s}$ is $A$-disjoint, then so is $\pi$.

**Lemma 3.6.** *For $A \in \mathbb{R}^{k \times n}$ and $p \geq 2$, the set of $A$-disjoint $p$-partitions equals the set of $p$-partitions associated with lists of $\binom{p}{2}$ $A$-disjoint 2-partitions.*

*Proof.* Let $\pi = (\pi_1, \ldots, \pi_p)$ be an $A$-disjoint $p$-partition. Consider any pair $1 \leq r < s \leq p$. Since $\mathrm{conv}(A^{\pi_r})$ and $\mathrm{conv}(A^{\pi_s})$ are disjoint, there exists an oriented hyperplane $H_{r,s}$ which contains no column of $A$ and satisfies $A^{\pi_r} \subset H_{r,s}^<$ and $A^{\pi_s} \subset H_{r,s}^>$. Let $\pi^{r,s} := (\pi_1^{r,s}, \pi_2^{r,s})$ be the $A$-disjoint 2-partition defined by $\pi_1^{r,s} := \{i \in [n] : A^i \in H_{r,s}^<\}$ and $\pi_2^{r,s} := \{i \in [n] : A^i \in H_{r,s}^>\}$. Let $\pi'$ be the $p$-tuple associated with the $\pi^{r,s}$. Then the $\pi'$ are pairwise disjoint and for $i = 1, \ldots, p$ we have
$$\pi_i \subseteq \left( \cap_{j=i+1}^{p} \pi_1^{i,j} \right) \bigcap \left( \cap_{j=1}^{i-1} \pi_2^{j,i} \right) = \pi_i'.$$
Since $[n] = \cup_{i=1}^p \pi_i \subseteq \cup_{i=1}^p \pi_i'$, it follows that $\pi = \pi'$ is the $p$-tuple associated with the constructed list of $\binom{p}{2}$ $A$-disjoint 2-partitions. □

Lemmas 3.5 and 3.6 imply our final preparatory lemma.

**Lemma 3.7.** *For $A \in \mathbb{R}^{k \times n}$ and $p \geq 2$, the set $\Pi_A^p$ of $A$-generic $p$-partitions equals the set of $p$-partitions associated with lists of $\binom{p}{2}$ $A$-generic 2-partitions.*

*Proof.* By Lemma 3.5 there is an $\epsilon$ such that $\Pi_A^2$ equals the set of $A(\epsilon)$-disjoint 2-partitions and $\Pi_A^p$ equals the set of $A(\epsilon)$-disjoint $p$-partitions. The lemma follows by applying Lemma 3.6 to $A(\epsilon)$. □

## 4. Vertex Enumeration and Optimization

We now use the facts established in the previous section to provide an upper bound on the number of vertices of any Shaped Partition Polytope. We then proceed to give efficient algorithms for the enumeration of these vertices and for the solution of the Shaped Partition Problem.

**Lemma 4.1.** *For any $A \in \mathbb{R}^{k \times n}$ and positive integer $p$, the number of $A$-generic $p$-partitions satisfies $|\Pi_A^p| = O(n^{k\binom{p}{2}})$.*

*Proof.* If $p = 1$ the bound is trivial so assume $p \geq 2$. By Lemma 3.5 and Lemma 3.3 applied to the perturbation $A(\epsilon)$ with $\epsilon$ sufficiently small, if $n \leq k$ then $\Pi_A^2$ consists of all $2^n = O(2^k) = O(n^k)$ 2-partitions, whereas if $n > k$ then it contains at most two associated partitions per triple $(I, J^<, J^>)$ with $I$ a $k$-subset of $[n]$ and $(J^<, J^>)$ a 2-partition of $I$. So $|\Pi_A^2|$ is at most twice the number $2^k \binom{n}{k}$ of such triples. By Lemma 3.7, the set $\Pi_A^p$ contains at most one partition per each list of $\binom{p}{2}$ partitions from $\Pi_A^2$, hence, as claimed,
$$|\Pi_A^p| \leq |\Pi_A^2|^{\binom{p}{2}} \leq \left( 2^{k+1} \binom{n}{k} \right)^{\binom{p}{2}} = O(n^{k\binom{p}{2}}). \quad \square$$

For $A \in \mathbb{R}^{k \times n}$ and nonempty set $\Lambda$ of $p$-shapes let $\Pi_A^\Lambda$ be the subset of $\Pi_A^p$ of $A$-generic $\Lambda$-partitions. For $A \in \mathbb{R}^{k \times n}$ let $\underline{A}$ be obtained as before by appending a new row $[1, \ldots, n]$ to $A$.

**Lemma 4.2.** *Let $A \in \mathbb{R}^{k \times n}$ and let $\Lambda$ be a nonempty set of $p$-shapes of $n$. Any vertex of $\mathcal{P}_A^\Lambda$ equals the $A$-matrix $A^\pi$ of some partition $\pi$ in $\Pi_{\underline{A}}^\Lambda$.*

*Proof.* Consider any vertex $A^\pi$ of $\mathcal{P}_A^\Lambda$. By Lemma 3.1, the $\Lambda$-partition $\pi$ is $A$-separable. By Lemma 3.2, there is another $\Lambda$-partition $\underline{\pi}$ such that $A^\pi = A^{\underline{\pi}}$ and $\underline{\pi}$ is $\underline{A}$-disjoint. By Lemma 3.5, $\underline{\pi} \in \Pi_{\underline{A}}^p$ hence is in $\Pi_{\underline{A}}^\Lambda$. □

We can now conclude an upper bound on the number of vertices of any Shaped Partition Polytope, which suffices for our purposes here. In [9] we use a more delicate argument which is based on treating the multisets of columns of $A$ in the original space, without lifting to $(k+1)$-space, and derive the more accurate bound of $O(n^{k\binom{p}{2}})$ on the number of vertices.

**Theorem 4.3.** *For any $A \in \mathbb{R}^{k \times n}$, positive integer $p$ and nonempty set $\Lambda$ of $p$-shapes of $n$, the number of vertices of the Shaped Partition Polytope $\mathcal{P}_A^\Lambda$ is $O(n^{(k+1)\binom{p}{2}})$.*

*Proof.* Let $\underline{A}$ be as above. By Lemmas 4.1 and 4.2 above,
$$|\text{vert}(\mathcal{P}_A^\Lambda)| \leq |\Pi_{\underline{A}}^\Lambda| \leq |\Pi_{\underline{A}}^p| = O(n^{(k+1)\binom{p}{2}}). \quad \square$$

Note that this theorem stands in contrast with the fact that typically the number of $\Lambda$-partitions is exponential in $n$ even for fixed $p$.

We next describe an algorithm for the enumeration of $\Pi_A^p$.

**Lemma 4.4.** *There is an algorithm that, given positive integers $p, k, n$ and a matrix $A \in \mathbb{R}^{k \times n}$, enumerates the set $\Pi_A^p$ of $A$-generic $p$-partitions using $O(p^2 \cdot n^{k\binom{p}{2}+1})$ real arithmetic operations. For rational $A$ and fixed $p$ and $k$ this is a strongly polynomial time algorithm.*

*Proof.* If $p = 1$ then let $\Pi_A^p := \{([n])\}$ consists of the single $p$-partition $([n])$. Suppose now $p \geq 2$. First, we describe how to produce the set $\Pi_A^2$ of $A$-generic 2-partitions. For each $\epsilon > 0$ let $A(\epsilon) \in \mathbb{R}^{k \times n}$ be a perturbation of $A$ obtained, as before, by setting $A(\epsilon)^i := A(\epsilon)^i + \epsilon M_k^i$, and let $\bar{A}(\epsilon)^i$ be obtained, as before, by appending a first coordinate 1 to $A(\epsilon)^i$. If $n \leq k$ then, by Lemma 3.3, $\Pi_A^2$ consists of all $2^n = O(2^k) = O(n^k)$ 2-partitions. Suppose now $n > k$. Let $I = \{i_1, \ldots, i_k\}$ be any $k$-subset of $[n]$ with $i_1 < \cdots < i_k$. For each $i \notin I$ let
$$D^i(\epsilon) := \det\left[\bar{A}(\epsilon)^{i_1}, \ldots, \bar{A}(\epsilon)^{i_k}, \bar{A}(\epsilon)^i\right] = \sum_{j=0}^k D_j^i \epsilon^j.$$

As explained in the proof of Lemma 3.5, $D^i(\epsilon)$ is a nonzero polynomial of degree $k$ in $\epsilon$ whose sign $D^i(\epsilon)$ is the same *generic* sign $\text{sign}(D^i)$ for all sufficiently small $\epsilon > 0$ and equals the sign of the first nonzero coefficient among $D_0^i, \ldots, D_k^i$. To compute these coefficients, substitute $\epsilon = 0, \ldots, k$ and compute each of the numerical determinants $D^i(0), \ldots, D^i(k)$. Then solve the system of linear equations $\sum_{j=0}^{k} \epsilon^j D_j^i = D^i(\epsilon)$, $\epsilon = 0, \ldots, k$ in the indeterminates $D_0^i, \ldots, D_k^i$. Since the defining matrix of this system is a Van der Monde matrix hence nonsingular, $D_0^i, \ldots, D_k^i$ are uniquely determined and so is $\text{sign}(D^i)$. Each of the numerical determinants $D^i(\epsilon)$ can be computed, and the system of equations solved, by Gaussian Elimination which takes $O(k^3)$ arithmetic operations. So each $\text{sign}(D^i)$ is computable in $O(k^4)$ operations. Therefore $I^< = \{i \in [n] : \text{sign}(D^i) < 0\}$ and $I^> = \{i \in [n] : \text{sign}(D^i) > 0\}$ are computable in $O(nk^4)$ operations, as are the two 2-partitions $\pi^<$ and $\pi^>$ associated with each triple $(I, J^<, J^>)$. As there are $2^k \binom{n}{k}$ such triples, the set $\Pi_A^2$ of $A$-generic 2-partitions can be computed in $2^k \binom{n}{k} nk^4 = O(n^{k+1})$ operations.

To construct $\Pi_A^p$, produce, using $O(p^2 n)$ operations, from each list of $\binom{p}{2}$ partitions from $\Pi_A^2$, the associated $p$-tuple $\pi$. Then test if it is a partition (i.e. if $\cup_{i=1}^{p} \pi_i = [n]$). As there are at most $(2^k \binom{n}{k})^{\binom{p}{2}}$ such lists, the total number of arithmetic operations is bounded by $p^2 n (2^k \binom{n}{k})^{\binom{p}{2}} = O(p^2 \cdot n^{k\binom{p}{2}+1})$ which subsumes the work for constructing $\Pi_A^2$ and is the claimed bound. Lemma 3.7 guarantees that each $A$-generic $p$-partition of shape $\Lambda$ is included. The claim about rational data follows from the fact that Gaussian Elimination admits a strongly polynomial time algorithm (see e.g. [16]). □

We can now provide the algorithm for the vertex enumeration of any Shaped Partition Polytope. It suffices for the algorithm that the set of shapes $\Lambda$ (whose size may be exponential in the rest of the data) is presented by an oracle that on query shape $\lambda$ answers whether or not $\lambda \in \Lambda$.

**Theorem 4.5.** *There is an algorithm that, given positive integers $p, k, n$, matrix $A \in \mathbb{R}^{k \times n}$, and an oracle presentation of a nonempty set $\Lambda$ of $p$-shapes of $n$, enumerates the set of vertices of the Shaped Partition Polytope $\mathcal{P}_A^\Lambda$ using $O(n^{k^2 p^3})$ real arithmetic operations and oracle queries. For rational $A$ and fixed $p$ and $k$ this is a strongly polynomial time oracle algorithm.*

*Proof.* Let $\underline{A} \in \mathbb{R}^{(k+1) \times n}$ be as before. Use the algorithm of Lemma 4.4 applied to $\underline{A}$ to construct the set of $\underline{A}$-generic $p$-partitions $\Pi_{\underline{A}}^p$. For each partition $\pi$ in this set, test if it is a $\Lambda$-partition by querying the $\Lambda$-oracle about its shape $|\pi|$. In this way, produce the set $\Pi_{\underline{A}}^\Lambda$ of $\underline{A}$-generic $\Lambda$-partitions. Then construct the set of $k \times p$ matrices $U := \{A^\pi : \pi \in \Pi_{\underline{A}}^\Lambda\}$ with multiple copies identified. This set $U$ is contained in $\mathcal{P}_A^\Lambda$ and by Lemma 4.2 contains the set of vertices of $\mathcal{P}_A^\Lambda$. So $u \in U$ will be a vertex precisely when it is not a convex combination of other elements of $U$. Testing it is a Linear Programming problem. To solve it, we shall appeal to Carathéodory's theorem and exhaustive search rather than invoke an LP-algorithm.

This allows to deal with real data and, for rational data and fixed $k, p$, obtain a *strongly* polynomial time procedure. So to test if $u \in U$ is a vertex of $\mathcal{P}_A^\Lambda$ do the following. Let $d := \dim(\text{aff}(U \setminus \{u\})) \leq kp + 1$. For each of the $\binom{|U|-1}{d}$ subsets $\{u_1, \ldots, u_d\}$ of $U \setminus \{u\}$, check if it is an affine basis of $\text{aff}(U \setminus \{u\})$ by checking if $\dim(\bar{u}_1, \ldots, \bar{u}_d) = d$; if it is, obtain the unique expression $\bar{u} = \sum_{i=1}^d \mu_i \bar{u}_i$ of $u$ as affine combination of the $u_i$. If $\mu_1, \ldots, \mu_d \geq 0$ then $u$ is in fact a convex combination of the $u_i$ hence not a vertex of $\mathcal{P}_A^\Lambda$. So $u$ is a vertex of $\mathcal{P}_A^\Lambda$ if and only if for each affine basis we get some $\mu_i < 0$. To compute $d$, to compute $\dim(\bar{u}_1, \ldots, \bar{u}_d)$, and to solve for the $\mu_i$, use Gaussian Elimination again.

By Lemma 4.4 applied to $\underline{A}$, the set $\Pi_{\underline{A}}^p$ can be constructed in $O(n^{(k+1)p^2})$ operations. The number of oracle queries is $|\Pi_{\underline{A}}^p| = O(n^{(k+1)\binom{p}{2}})$, as guaranteed by Lemma 4.1, and the number of matrices in $U$ satisfies $|U| \leq |\Pi_{\underline{A}}^\Lambda| \leq |\Pi_{\underline{A}}^p| = O(n^{(k+1)\binom{p}{2}})$. As the effort of constructing $\Pi_{\underline{A}}^\Lambda$ is absorbed in the vertex testing procedure, we have that the total number of arithmetic operations and oracle queries is bounded by $|U|\binom{|U|-1}{kp+1}(kp+1)^3 = O(n^{k^2 p^3})$ as claimed. The claim for rational data follows from the strongly polynomial time computability of Gaussian Elimination. □

Finally, we provide an algorithmic solution of the Shaped Partition Problem we started with. It suffices that the convex functional $C$ is presented by an oracle that on query $A^\pi$ with $\pi$ a $\Lambda$-partition answers with $C(A^\pi)$. The oracle for $C$ will be called *M-guaranteed* if $C(A^\pi)$ is guaranteed to be a rational number whose binary encoding is no longer than $M$ bits for any $\Lambda$-partition $\pi$.

**Theorem 4.6.** *There is an algorithm that, given positive integers $p, k, n$, matrix $A \in \mathbb{R}^{k \times n}$, an oracle presentation of a nonempty set $\Lambda$ of $p$-shapes of $n$, and an oracle presentation of a convex functional $C$ on $\mathbb{R}^{k \times p}$, solves the Shaped Partition Problem of finding a $\Lambda$-partition $\pi^*$ which maximizes $f(\pi) = C(A^\pi)$ using $O(n^{kp^2})$ real arithmetic operations and oracle queries. For rational $A$, fixed $p$ and $k$, and $M$-guaranteed oracle for $C$, the running time of this oracle algorithm is strongly polynomial in $n, M$ and the bit size of $A$.*

*Proof.* Let $\underline{A}$ be as before. Use the algorithm of Lemma 4.4 applied to $\underline{A}$ to construct the set of partitions $\Pi_{\underline{A}}^p$, using $O(p^2 \cdot n^{(k+1)\binom{p}{2}+1})$ arithmetic operations. Then filter out its subset $\Pi_{\underline{A}}^\Lambda$ of $\underline{A}$-generic $\Lambda$-partitions by querying the $\Lambda$-oracle on each of the $|\Pi_{\underline{A}}^p| = O(n^{(k+1)\binom{p}{2}})$ partitions in $\Pi_{\underline{A}}^p$. Since $C$ is convex, there will be an optimal $\Lambda$-partition $\pi$ for the Shaped Partition Problem for which $A^\pi$ is a vertex of $\mathcal{P}_A^\Lambda$, hence, by Lemma 4.2, an optimal $\pi \in \Pi_{\underline{A}}^\Lambda$. For each $\pi \in \Pi_{\underline{A}}^\Lambda$ compute $A^\pi$ (using $O(kn)$ operations) and query the $C$-oracle for the value $C(A^\pi)$. Pick $\pi^*$ to be the one attaining maximum value. The number of operations and queries to the $\Lambda$-oracle and $C$-oracle is bounded, as claimed, by $O(kn \cdot n^{(k+1)\binom{p}{2}}) + O(p^2 \cdot n^{(k+1)\binom{p}{2}+1}) = O(n^{kp^2})$.

A POLYTIME ALGORITHM FOR SHAPED PARTITION PROBLEMS AND POLYTOPES    13For rational data, since the $C$-oracle is $M$-guaranteed, each of the comparisons of its answers needed to determine the optimal $\pi^*$ is doable in time linear in $M$. Since $\Pi_A^\Lambda$ is strongly polynomial time computable for fixed $p,k$, the claim about strongly polynomial time solution of the Shaped Partition Problem follows.  □

## Acknowledgment

Shmuel Onn thanks the Mathematical Sciences Research Institute at Berkeley for its support while part of this research was done.Shmuel Onn thanks the Mathematical Sciences Research Institute at Berkeley for its support while part of this research was done.

Frank K. Hwang, Department of Applied Mathematics, Chaiotung University, Hsinshu, Taiwan.

Shmuel Onn, Department of Operations Research, Faculty of Industrial Engineering and Management, Technion - Israel Institute of Technology, 32000 Haifa, Israel.

Uriel G. Rothblum, Department of Operations Research, Faculty of Industrial Engineering and Management, Technion - Israel Institute of Technology, 32000 Haifa, Israel.